\documentclass[11pt]{article}
\usepackage{amsmath, amssymb, amsfonts, enumerate}
\newtheorem{theorem}{Theorem}[section]

\newtheorem{proposition}[theorem]{Proposition}
\newtheorem{corollary}[theorem]{Corollary}
\newtheorem{lemma}[theorem]{Lemma}

\newtheorem{problem}[theorem]{Problem}
\newtheorem{remark}[theorem]{Remark} 

\def\IR{{\rm I\!R}} 
\def\IC{\hbox{\rm C\kern-.43em
       \vrule depth 0ex height 1.4ex width .05em\kern.41em}}
\def\IQ{\hbox{\rm Q\kern-.43em
       \vrule depth 0ex height 1.4ex width .05em\kern.41em}}

\def\bC{{\bf C}}   
 \def\bR{{\bf R}}

\def\bR{{\bf R}}

\def\qed{\hfill\vbox{\hrule width 6 pt
\hbox{\vrule height 6 pt width 6 pt}} \medskip}

\begin{document}
\openup .8\jot
\title{On the boundary of weighted numerical ranges}
\author{Wai-Shun Cheung}
\date{}
\maketitle

\begin{abstract}

In this article, we introduce the weighted numerical range which is a unified approach to study the $c$-numerical range and the rank $k$ numerical range.  If the boundaries of weighted numerical ranges of two matrices (possibly of different sizes) overlap at sufficiently many points, then the two matrices share common generalized eigenvalues. 

\end{abstract}



\section{Introduction}

Let $M_n$ denote the space of all $n\times n$ complex matrices and $\IR^n$ the set of all real $n$-tuples.  For any $A\in M_n$, we denote $\lambda_1(A), \ldots, \lambda_n(A)$ the $n$ eigenvalues of $A$.  In the case that $A$ is hermitian, we assume that $\lambda_1(A)\ge \lambda_2(A)\ge\cdots\ge\lambda_n(A)$. We also define $\displaystyle H_\theta(A)=\frac{e^{i\theta}A+e^{-i\theta}A^*}{2}$ for $\theta\in [0,2\pi)$.

Numerical range is an extensively studied subject area.  The study starts with the classical numerical range of $A\in M_n$ which is defined as
$$W(A)=\{x^*Ax\ :\ x\in \IC^n, x^*x=1\}$$
which is a compact set containing all the eigenvalues of $A$.  $W(A)$ is a convex set by the famous Toeplitz-Housedorf Theorem \cite{T,H}. A nice discussion can be found in \cite[Chapter 1]{HJ}.   There are many papers related to the boundary of the classical numerical range \cite{K, RS, GW1, GW2, Wu, CL}.

There are many different generalizations of the classical numerical range.

For $1\le k\le n$, the $k$-numerical range of $A$, introduced by Halmos \cite{Hal} is defined as
$$W_k(A)=\left\{\sum_{j=1}^k\frac{1}{k}x_j^*Ax_j\ :\ x_1,\ldots,x_k\in \IC^n \mbox{ are orthonormal}\right\}$$
which is proved to be convex by Berger \cite{B}.  Note that $W_1(A)=W(A)$. Indeed $W_k(A)\subseteq W(A)$.

For any $c=(c_1,c_2,\ldots,c_n)^t\in \IR^n$, the $c$-numerical range of $A$, first introduced by Marcus \cite{M}, is defined as
$$W_c(A)=\left\{\sum_{j=1}^nc_jx_j^*Ax_j\ :\ x_1,\ldots,x_n\in \IC^n \mbox{ are orthonormal}\right\}.$$
Westwick \cite{W} proved it is a convex set.  Indeed, if $c_1=\cdots=c_k=1/k,c_{k+1}=\cdots=c_n=0$ then $W_c(A)=W_k(A)$. Therefore the $c$-numerical range is a natural generalization of the $k$-numerical range.  A survey can be found in \cite{Li}.  There are papers related to the boundary of the $k$-numerical range and the $c$-numerical range \cite{CL, MMF, LST}.

A seemingly different generalization is the rank $k$ numerical range.

For any $1\le k\le n$, the rank $k$ numerical range of $A$, first introduced by  Choi, Kribs and \.{Z}yczkowsk \cite{CKZ, CKZ1}, is defined as
$$\Lambda_k(A)=\{\lambda\in\IC\ :\ PAP=\lambda P \mbox{ for some rank-}k \mbox{ orthogonal projection }P\}.$$
Note that $\Lambda_1(A)=W(A)$. The rank $k$ numerical of $A$ was proven to be convex by Woerdeman \cite{Wo} and Li and Sze \cite{LS} independently. The rank $k$-numerical range is a relatively new generalized numerical range and it is a hot topic partly due to its connection to quantum computing, e.g. \cite{LP, Many}.  There are papers on the boundary of the rank $k$ numerical range \cite{CN, GLPS, LPS, C}, but not much geometric properties are known.

As it turns out, all the three generalizations have similar reformulations.  In this article, we are going to introduce the weighted numerical range, which is a unified approach to the $c$-numerical range and the rank $k$ numerical range.  We will prove a theorem that if the boundaries of the weighted numerical ranges of two matrices have many intersection points, then the two matrices have common generalized eigenvalues.  The size of the two matrices could be different. Applying the theorem and using known results on classical numerical range, we can deduce some  properties about the $c$-numerical range and the rank $k$ numerical range.

\section{Weighted Numerical Ranges}
From \cite[Chapter 1]{HJ}, we have the following equivalent expression for $W(A)$:
$$W(A)=\bigcap_\theta \left\{v\in \IC\ : \mbox{Re }e^{i\theta}v \le \lambda_1(H_\theta(A))\right\}.$$

There is a similar equivalent expression for $W_c(A)$:
$$W_c(A)=\bigcap_\theta \left\{v\in \IC\ : \mbox{Re } e^{i\theta}v \le c_{\sigma(1)}\lambda_1(H_\theta(A))+\cdots+c_{\sigma(n)}\lambda_n(H_\theta(A))\right\}$$
where $\sigma\in S_n$ such that $c_{\sigma(1)}\ge \ldots\ge c_{\sigma(n)}$. 

Li and Sze \cite{LS} proved the convexity of $\Lambda_k(A)$ by showing that there is also a similar equivalent expression of $\Lambda_k(A)$:
$$\Lambda_k(A)=\bigcap_\theta \left\{v\in \IC\ : \mbox{Re }e^{i\theta}v \le \lambda_k(H_\theta(A))\right\}.$$

Inspired by the alternative expressions, we define a new type of generalized numerical range.  
For any $c=(c_1,c_2,\ldots,c_n)^t\in \IR^n$, we define
\begin{equation}
W(A;c)=\bigcap_\theta \left\{v\in \IC\ : \mbox{Re } e^{i\theta}v \le c_1\lambda_1(H_\theta(A))+\cdots+c_n\lambda_n(H_\theta(A))\right\}.
\end{equation}
Follow the suggestion of Karol  \.{Z}yczkowski, we will call it the {\it weighted numerical range}.

Lets start with some simple properties of the weighted numerical ranges.

\begin{proposition}
Let $c\in\IR^n$ and $A\in M_n$. The weighted numerical range of $A$ has the following properties:
\begin{enumerate}
\item If $W(A;c)$ is nonempty, then it is a complex convex set.
\item When $c=e_k=(0,\ldots,0,1,0,\ldots,0)$, $W(A;c)=\Lambda_k(A)$. 
\item $W_c(A)=W(A; (c_{\sigma(1)}, \ldots, c_{\sigma(n)}))$ where $\sigma$ is a permutation such that  $c_{\sigma(1)}\ge \ldots\ge c_{\sigma(n)}$.
\item If $A$ is normal, then $W(A;c)$ is a polygonal disc.
\item If $A$ is hermitian, then $W(A;c)$ is the real segment $\{x\ :\ \sum_{j=1}^n c_{n+1-j}\lambda_j(A)\le x\le \sum_{j=1}^n c_j\lambda(A)\}$ which can be empty.
\item $W(\gamma A+\lambda I; c)=\gamma W(A;c)+\lambda(\sum_{j=1}^nc_j)$.
\end{enumerate}
\end{proposition}

{\it Proof.} (1), (2), (3) and (6) follow directly from the definiton.

Suppose $A$ is normal with eigenvalues $u_j+iv_j$ for $j=1,\ldots, n$.  The eigenvalues of $H_\theta(A)$ are therefore $u_j\cos\theta-v_j\sin\theta$, $j=1,\ldots, n$.  We define $\sigma=\sigma(\theta)$ to be a permutation such that $u_{\sigma(j)}\cos\theta-v_{\sigma(j)}\sin\theta$, $j=1,\ldots,n$ are in decreasing order.
The equality
$$\mbox{Re }e^{i\theta}(x+i y) \le c_1\lambda_1(H_\theta(A))+\cdots+c_n\lambda_n(H_\theta(A))$$
is therefore equivalent to
$$x \cos\theta - y\sin\theta \le \sum_{j=1}^n c_{\sigma^{-1}(j)} (u_j\cos\theta-v_j\sin\theta).$$
If $\cos\theta>0$, it becomes
$$x- \sum_{j=1}^n c_{\sigma^{-1}(j)} \le (y-\sum_{j=1}^n c_{\sigma(j)} v_j)\tan\theta.$$
If $\cos\theta<0$, it becomes
$$x- \sum_{j=1}^n c_{\sigma^{-1}(j)} \ge (y-\sum_{j=1}^n c_{\sigma(j)} v_j)\tan\theta.$$
If $\theta=\pi/2$, it becomes
$$y \le \sum_{j=1}^n c_{\sigma^{-1}(j)} v_j.$$
If $\theta=-\pi/2$, it becomes
$$y \ge \sum_{j=1}^n c_{\sigma^{-1}(j)} v_j.$$
Therefore for each $\sigma$, the set of corresponding inequalities can be reduced to a set of at most four inequalities. That is, for $x+i y\in W(A;c)$, $x$ and $y$ need to satisfies atmost $4n!$ inequalities.  In other words, $W(A;c)$ is a polygonal disc with at most $4n!$ sides. Hence we have (4).

Suppose $A$ is hermitian. Apply the proof of (4) and note that only two permutations are relevant, we have (5).   
\qed

\section{The $c$-values and $c$-polynomial}
Let $c=(c_1,\ldots,c_n)^t\in \IR^n$ and $A\in M_n$.  Suppose $c_{i_1}, \ldots, c_{i_r}$ are all the nonzero entries. We define
$$\lambda_{j_1,\ldots,j_r}(A;c)=c_{i_1}\lambda_{j_1}(A)+\cdots+c_{i_r}\lambda_{j_r}(A)$$
and call it a $c$-value of $A$.
The $c$-polynomial of $A$,  $p(A;c)(t)$, is the polynomial which takes all generic distinct $c$-values as root.

For example, if $c=(1,0,1,2)$ and that $a,b,c,d$ are the eigenvalues of $A$, then $a+b+2c, a+b+2d, a+c+2b, a+c+2d, a+d+2b, a+d+2c, b+c+2a, b+c+2d, b+d+2a, b+d+2c, c+d+2a, c+d+2b$ are all the $c$-values of $A$ and $p(A;c)(t)$ is a polynomial of degree 12.

We introduce two more notations.  Let $r(A;c)(x,y,t)=p(xA+yA^*;c)(t)$ which is  a polynomial in $x,y,t$ and that $r(A;c)(1,0,t)=p(A;c)(t)$. We also write $\deg(A;c)=\deg p(A;c)(t)$. 

Note that if $r=\deg (A;c)$ then the coefficients of $t^k$ in $r(A;c)=p(xA+yA^*;c)(t)$ are symmetric homogeneous functions of degree $r-k$ on the $c$-values of $xA+yA^*$, and hence symmetric homogeneous functions of degree $r-k$ on the coefficients of the characteristic polynomial of $xA+yA^*$, and hence symmetric homogeneous functions of degree $r-k$ on the entries of $xA+yA^*$. Therefore $r(A;c)$ is a homogeneous function of degree $r$ and that $\deg(A;c)=
\deg r(A;c)(x,y,t)$.

Apply Bezout's Theorem, we have two lemmas.
\begin{lemma} \label{lemma1}
Let $A\in M_n$, $B\in M_m$, $c\in\IR^n$ and $d\in\IR^m$.  If $r(A;c)$ and $r(B;d)$ have more than $\deg(A;c) \deg(B;d)$ common roots in the projective plane, then $r(A;c)$ and $r(B;d)$ have a common factor.  Consequently there exists a $c$-value of $A$ which is also a $d$-value of $B$.
\end{lemma}

\begin{lemma} \label{lemma2}
Let $A\in M_n$, $B\in M_m$, $c\in\IR^n$ and $d\in\IR^m$.  If $r(A;c)$ and $r(B;d)$ have more than $\deg(A;c) \deg(B;d)$ common roots in the projective plane and that $r(A;c)$ is irreducible, then $r(A;c)$ is a factor of $r(A;c)$.  Consequently all the $c$-values of $A$ are $d$-values of $B$.
\end{lemma}

\section{Supporting Lines and Common Boundary Points}

We need two technical lemmas for further discussion.  The first lemma is trivial, it tells us when a supporting line is of a special form.

\begin{lemma} \label{lemma3}
A supporting line of $W(A;c)$ is of the form 
$$\left\{v\in \IC\ : \mbox{Re }e^{i\theta}v = c_1\lambda_1(H_\theta(A))+\cdots+c_n\lambda_n(H_\theta(A))\right\}$$
whenever
\begin{enumerate}
\item[(a)] the supporting line is tangent to the boundary of $W(A;c)$ at a differentiable point; or
\item[(b)] $c_1\ge c_2\ge \cdots\ge c_n$. In other words, every supporting line of $W_c(A)$ is of this special form.
\end{enumerate}
\end{lemma}

The second lemma states that if there are three common boundary points of two weighted numerical ranges, then there must be a common supporting line.

\begin{lemma} \label{lemma4}
Let $A\in M_n$, $B\in M_m$, $c\in\IR^n$ and $d\in\IR^m$.
Suppose $z_1, z_2, z_3\in \partial W(A;c)\cap  \partial W(B;d)$ are positioned along the boundary of $W(A;c)$ (or $W(B;d)$) in an anticlockwise way.  Let $\omega_l=Arg(i\overline{z_{l+1}-z_{l}})$, $l=1,2$. then there exists $\phi\in [\omega_2,\omega_1]$ (defined as $[\omega_2, 2\pi]\cup [0, \omega_1]$ if $\omega_1<\omega_2]$) such that 
$$\sum_{j=1}^n c_j \lambda_j(H_\phi(A))=\sum_{k=1}^m d_k \lambda_k(H_\phi(A)).$$
In other words, there is a common supporting line in the direction which is between the perpendicular bisectors of $[z_1,z_2]$ and $[z_2,z_3]$.
Furthermore $\phi$ can be chosen so that $\phi\in (\omega_1,\omega_2)$ or at least one of $[z_1,z_2]$ and $[z_2,z_3]$ is part of $\partial W(A;c)\cap \partial W(B;d)$.
\end{lemma}

{\it Proof.} Let $L:\{z\ :\  \mbox{Re } e^{i\theta}z=\sum_{k=1}^m d_k \lambda_k(H_\theta(B))\}$ be a supporting line of $W(B;d)$ passing $z_2$.

Let $z_2-z_1=r e^{i(\pi/2-\omega_1)}$, then
$$\mbox{Re }e^{i\theta}(z_2-z_1)=r\cos(v+\pi/2-\omega_1)=r\sin(\omega_1-v)$$
is nonnegative if: when $0\le\omega_1\le \pi$, then  $0\le \theta\le \omega_1$ or $\pi+\omega_1\le \theta\le 2\pi$; when $\pi<\omega_1\le 2\pi$, then $\omega_1-\pi \le \theta\le \omega_1$.

Likewise $\mbox{Re }e^{i\theta}(z_2-z_3)$ is nonnegative if: when $0\le\omega_2\le\pi$, then $\omega_2\le \theta\le \pi+\omega_2$; when $\pi<\omega_2\le 2\pi$, then $0\le \theta\le \omega_2-\pi$ or $\omega_2\le\theta\le 2\pi$. 

If $\omega_2\le \omega_1$ then $\omega_1\le\theta\le \omega_2$. 
If $\omega_2>\omega_1$ and $\theta>\omega_1$ then $\theta>\omega_1+\pi$ and hence $\theta>\omega_2$.
Thus $\theta\in [\omega_2, \omega_1]$. 

Now $\sum_{j=1}^n c_j \lambda_j(H_\theta(A))\ge \sum_{k=1}^m d_k \lambda_k(H_\theta(B))$ or otherwise $z_2\notin W(A;c)$.

Likewise, we can find $\varphi\in [\omega_2,\omega_1]$ such that $\sum_{j=1}^n c_j \lambda_j(H_\varphi(A))\le \sum_{k=1}^m d_k \lambda_k(H_\varphi(B))$.

By the continuity of $\lambda_j(H_\theta(A))$ and $\lambda_k(H_\theta(B))$, there exists $\phi$ between $\theta$ and $\varphi$ such that $\sum_{j=1}^n c_j \lambda_j(H_\phi(A))=\sum_{k=1}^m d_k \lambda_k(H_\phi(A)).$

If $\phi=\omega_1$, then either $\phi=\omega_1=\theta$ or $\phi=\omega_1=\varphi$. In both case, it implies that $[z_1,z_2]\subseteq \partial W(A;c)\cap \partial W(B;d)$.

Similarly for the case $\phi=\omega_2$.
\qed

\section{Main Theorems}

We are ready to state the main results.

\begin{theorem} \label{main}
Let $A\in M_n$, $B\in M_m$, $c\in\IR^n$ and $d\in\IR^m$. Suppose
$$c_1\lambda_1(H_\theta(A))+\cdots+c_n\lambda_n(H_\theta(A))=
d_1\lambda_1(H_\theta(B))+\cdots+d_m\lambda_m(H_\theta(B))$$
for $\deg(A;c) \deg(B;d)+1$ $\theta$'s, then there exists a $c$-value of $A$ which is also a $d$-value of $B$.  Furthermore, if $r(A;c)$ is irreducible, then all the $c$-values of $A$ are $d$-values of $B$.
\end{theorem}

{\it Proof. }The condition
$$r=c_1\lambda_1(H_\theta(A))+\cdots+c_n\lambda_n(H_\theta(A))
=d_1\lambda_1(H_\theta(B))+\cdots+d_m\lambda_m(H_\theta(B))$$
implies that a $c$-value of $H_\theta(A)$ is a $d$-value of $H_\theta(B)$,
which in turns implies that $(e^{i\theta}, e^{-i\theta}, r)$ is a common root of $r(A;c)$ and $r(B;d)$.
The result then follows Lemma~\ref{lemma1} and Lemma~\ref{lemma2}.
\qed

Apply Lemma~\ref{lemma3} and Theorem~\ref{main}, we prove a theorem on common supporting lines.

\begin{theorem} \label{thm1}
Let $A\in M_n$, $B\in M_m$, $c\in\IR^n$ and $d\in\IR^m$. If $W(A;c)$ and $W(B;d)$ have more than  $\deg(A;c) \deg(B;d)$ common supporting lines and the following two conditions are satisfied:
\begin{enumerate}
\item $W(A;c)=W_c(A)$, or each supporting line touches $\partial W(A;c)$ at a differentiable point;
\item $W(B;d)=W_d(B)$, or each supporting line touches $\partial W(B;d)$ at a differentiable point,
\end{enumerate}
then there exists a $c$-value of $A$ which is also a $d$-value of $B$.  Furthermore, if $r(A;c)$ is irreducible, then all the $c$-values of $A$ are $d$-values of $B$.
\end{theorem}

Apply Lemma~\ref{lemma4} and Theorem~\ref{main}, we prove a theorem on common boundary points.

\begin{theorem} \label{thm2}
Let $A\in M_n$, $B\in M_m$, $c\in\IR^n$ and $d\in\IR^m$.  Suppose there are $z_1,\ldots,z_k \in \partial W(A;c)\cap \partial W(B;d)$ where $k=\deg(A;c) \deg(B;d)+1$ such that $[z_r,z_s]$ do not lie on $\partial W(A;c)\cap \partial W(B;d)$ for $r\ne s$, then there exists a $c$-value of $A$ which is also a $d$-value of $B$.  Furthermore, if $r(A;c)$ is irreducible, then all the $c$-values of $A$ are $d$-values of $B$.
\end{theorem}

{\it Proof.} Suppose the $k$ points are in an anticlockwise manner and that $z_{k+1}=z_1, z_{k+2}=z_2$.

By Lemma~\ref{lemma4}, $\{z_l, z_{l+1}, z_{l+2}\}$ defines an angle $\phi_l\in  (\omega_{l+1}, \omega_{l})$ where  $\omega_s=Arg(i\overline{z_{s+1}-z_{s}})$ such that
$$\sum_{j=1}^n c_j \lambda_j(H_\phi(A))=\sum_{k=1}^m d_k \lambda_k(H_\phi(A)).$$
Note that those $k$  $\phi$'s are distinct. Therefore, by Theorem~\ref{main}, the result follows. 
\qed

We have two simple consequences of Theorem~\ref{thm1}.

\begin{corollary} \label{thm3}
Let $A\in M_n$, $B\in M_m$, $c\in\IR^n$ and $d\in\IR^m$. Suppose there is a differentiable curve, which is not a straight line, lying on $\partial W(A;c) \cap \partial W(B;d)$, then there exists a $c$-value of $A$ which is also a $d$-value of $B$.  Furthermore, if $r(A;c)$ is irreducible, then all the $c$-values of $A$ are $d$-values of $B$.
\end{corollary}

\begin{corollary} \label{thm4}
Let $A\in M_n$, $B\in M_m$, $c\in\IR^n$ and $d\in\IR^m$.  If $W_c(A)=W_d(B)$ then there exists a $c$-value of $A$ which is also a $d$-value of $B$.
\end{corollary}

The next two corollaries relate to some old results \cite{ GW1, GW2, Wu, CL, C}.

\begin{corollary} \label{cor1}
Let $A\in M_n$ and $c\in \IR^n$. If $\partial W(A;c)$ contains $2\deg(A;c)+1$ points on a circle centered at $\alpha$, then $\alpha$ is a $c$-value of $A$ with multiplicity greater than $1$.

Consequently, if $\Lambda_k(A)$ contains a circular arc centered at $\alpha$, than $\alpha$ is an eigenvalue of $A$ with multiplicity greater then $1$.
\end{corollary}

{\it Proof.} Let $B=\begin{pmatrix}\alpha&2R\\0&\alpha\end{pmatrix}$ where $R$ is the radius of the arc and let $d=(1,0)$. Apply Theorem~\ref{thm2} and note that $r(B;d)$ is irreducible.
\qed 

\begin{corollary} \label{cor2}
Let $A\in M_n$ and $c\in \IR^n$. If $\partial W(A;c)$ contains $2\deg(A;c)+1$ points on an ellipse, then the two foci of the ellipse along the main axis are $c$-values of $A$.

Consequently, if $\Lambda_k(A)$ contains an elliptical arc,  then the two foci of the elliptical arc along the main axis are two eigenvalues of $A$.
\end{corollary}

{\it Proof.} Let $B=\begin{pmatrix}\alpha&R\\0&\beta\end{pmatrix}$ where $\alpha$ and $\beta$ are the foci of the ellipse and $R$ is a suitable number. Let $d=(1,0)$. Apply Theorem~\ref{thm2} and note that $r(B;d)$ is irreducible. 
\qed 

\begin{remark}
The bound $\deg(A;c)\deg(B;d)+1$ is sharp.  Let $A$ to be the $n\times n$ diagonal matrix with eigenvalues being the $n$ roots of unity and $B=\begin{pmatrix}0&2R\\0&0\end{pmatrix}$ where $R$ is slightly less than $1$, then $W(A)$ and $W(B)$ have exactly $2n$ common boundary points, but $A$ and $B$ have no common eigenvalues.
\end{remark}

We end this section with two known results with new proofs.

The first result is on the sharp point of $c$-numerical ranges.
\begin{corollary}
Let $A\in M_n$ and $c\in \IR^n$. If $W_c(A)$ has a sharp point $\alpha$ then $\alpha$ is a $c$-value of $A$.
\end{corollary}

{\it Proof.} Let $B=\alpha I_2$ and $d=(1,0)$. Apply Theorem~\ref{thm1}. \qed

The second result in \cite{Ta} relates to in \cite{LT} and a follow-up question listed in \cite[Section 9]{Li}.  

\begin{corollary}
Consider $A\in M_n(\bC)$. Suppose $W_c(A)$ is a circular disc centered at $0$ for any $c\in \bR^n$, then $A$ is nilpotent.
\end{corollary}
{\it Proof. } Let $a_1,\ldots,a_n$ be eigenvalues of $A$.  Suppose not all $a_j$'s are zero, then there exists $c=(c_1,\ldots,c_n)^t\in\bR^n$ such that $c_1a_{\sigma(1)}+\cdots+c_na_{\sigma(n)}\ne 0$ for all permutations $\sigma$.

Since $W_c(A)$ is a circular disc centered at $0$, we have $W_c(A)=W(\alpha E_{12})$ for some $\alpha\in\bR$. By Theorem~\ref{cor1}, there exists a permutation $\sigma$ such that
$c_1a_{\sigma(1)}+\cdots+c_na_{\sigma(n)}=0$.

We now have a contradiction.
\qed

\section{Open Questions}

\begin{problem}
What happens if there is a sharp point or a line segment on $\partial W(A;c)$?
\end{problem}

We know very little even for rank $k$ numerical range.

\begin{problem}
Could we get any meaningful results if $\partial W(A;c)\cap \partial W(B;d)$ contains a line segment?
\end{problem}

Again, we know very little for rank $k$ numerical range.

\begin{problem}
Suppose we know that $W(A;c)\subseteq W(B;d)$ and that $\partial W(A;c)\cap \partial W(B;d)$ contains sufficiently many points. Could we say more about the geometry of $W(A;c)$ and $W(B;d)$?
\end{problem}

Wu \cite{Wu} proved some nice results if $W(A)$ or $W(B)$ is a circular disc.  Cheung and Li \cite{CL} generalized Wu's results to elliptical disc and $k$-numerical range. Cheung \cite{C} obtained extended Wu's results to rank $k$ numerical range.  We believe that there should be some similar results for weighted numerical range.

\section*{Acknowledgement}
I would like to thank C.K. Li and Karol \.{Z}yczkowski for their valuable suggestions. I would also like to thank the referee for pointing out a mistake in the first version of Lemma~\ref{lemma4} and other errors.

I would also like to thank P.S. Lau for reminding me of \cite{Ta} after the paper is accepted by LAMA.

%C.K. Li has asked if I can arrive at something related to the sharp points on the boundary of the weighted numerical range. Sadly I am not yet to able to get any interesting result.

\medskip
\sc

Department of Mathematics, University of Hong Kong, Hong Kong.

\rm

E-mail addresses: cheungwaishun@gmail.com
\end{document}